\documentclass[12pt]{amsart}
\usepackage{amssymb}
\usepackage{times}
\usepackage{graphicx}

\theoremstyle{plain}

\numberwithin{equation}{section}

\newcommand{\calA}{\mathcal{A}}
\newcommand{\calB}{\mathcal{B}}

\newcommand{\calD}{\mathcal{D}}

\newcommand{\calH}{\mathcal{H}}

\newcommand{\calM}{\mathcal{M}}

\begin{document}
\title [The moduli space of 8 points on ${\bf P}^1$ ] {The moduli space of 8 points on ${\bf P}^1$ and automorphic forms}
\author{Shigeyuki Kond{$\bar{\rm o}$}}
\address{Graduate School of Mathematics, Nagoya University, Nagoya,
464-8602, Japan}
\email{kondo@math.nagoya-u.ac.jp}
\thanks{Research of the author is partially supported by
Grant-in-Aid for Scientific Research A-14204001, Japan}
\begin{abstract}
First we give a complex ball uniformization of the moduli space of 8 ordered points 
on the projective line by
using the theory of periods of $K3$ surfaces.  Next we give a projective model of this
moduli space by using automorphic forms on a bounded symmetric domain of type $IV$ 
which coincides with the one given by cross ratios of 8 ordered points of the projective line.
\end{abstract}
\maketitle

\rightline{\it Dedicated to Igor Dolgachev on  his 60th birthday}

\section{Introduction}
The main purpose of this paper is to give an application of the theory of automorphic forms on 
a bounded symmetric domain of type $IV$ due to Gritsenko and Borcherds \cite{Bor} for studying
the moduli spaces.  We consider the moduli space $P_1^8$ of semi-stable 8 ordered points on the
projective line.  It is known that $P_1^8$ is isomorphic to 
the Satake-Baily-Borel compactification of an arithmetic quotient of 5-dimensional complex
ball by using the theory of periods of a family of curves which are the 4-fold cyclic covers of
the projective line branched at eight points (\cite{DM}).  Recently
Matsumoto and Terasoma \cite{MT} gave an embedding of $P_1^8$ into ${\bf P}^{104}$ by
using the theta constants related to the above curves.  Their map coincides the one defined by 
the cross ratios of 8 points on the ${\bf P}^1$.  Here they used the fact that the complex
ball is canonically embedded in a Siegel upper half plane.  

In this paper, instead of the periods of curves, we use the periods of $K3$ surfaces.  In our case,
the complex ball is embedded in a bounded symmetric domain of type $IV$.
In fact, to each stable point from $P_1^8$, we associate a $K3$ surface with a non-symplectic
automorphism of order 4 (\S \ref{k3}).  
The period domain of these $K3$ surfaces is a 5-dimensional
complex ball $\calB$ (\S \ref{unif}).  This was essentially given in the paper \cite{Kon1}. 
By using this, we shall see that
$P_1^8$ is isomorphic to the Satake-Baily-Borel compactification ${\bar \calB}/\Gamma(1-i)$ of 
$ \calB/\Gamma(1-i)$ where $\Gamma(1-i)$ is an arithmetic subgroup of
a unitary group of a hermitian form of signature $(1,5)$ defined over the Gaussiann integers
(Theorems \ref{main}, \ref{Main}) .  The symmetry group $S_8$ of degree 8 naturally acts on $P_1^8$.  On the other hand, there exists an arithmetic subgroup $\Gamma$ acting on 
$\calB$ with $\Gamma/\Gamma(1-i) \cong S_8$.  The above isomorphism 
$P_1^8 \cong {\bar  \calB}/\Gamma(1-i)$ is $S_8$-equivariant.

Next we apply the theory of automorphic forms \cite{Bor} to this situation.  The main idea comes from the paper of Allcock and Freitag \cite{AF} in which they studied the same problem in the case of cubic surfaces.  
We shall show that there exists a 14-dimensional space of automorphic forms on $\calB$ which
gives an $S_8$-equivariant  map from the arithmetic quotient ${\bar \calB/\Gamma(1-i)}$ into 
${\bf P}^{13}$ (Theorem \ref{morphism}).
Under the identification $P_1^8 \cong  {\bar \calB/\Gamma(1-i)}$ 
we show that this map coincides with the one defined by the cross ratios of 8 points on the
projective line (Theorem \ref{cross}).  Thus our map coincides the one given by Matsumoto and 
Terasoma \cite{MT}.

In this paper, a {\it lattice} means a $\bf Z$-valued non-degenerate symmetric bilinear form on 
a free $\bf Z$-module of finite rank.  We denote by $U$ the even lattice defined by the matrix
$
\begin{pmatrix}0&1
\\1&0
\end{pmatrix}
$,
and by $A_m$, $D_n$ or $E_l$ the even negative definite lattice defined by the Dynkin matrix of
type $A_m$, $D_n$ or $E_l$ respectively.  If $L$ is a lattice and $m$ is an integer, 
we denote by $L(m)$ the lattice over the same $\bf Z$-module with the symmetric bilinear form
multiplied by $m$.  We also denote by $L^{\oplus m}$ the orthogonal direct sum of $m$ copies
of $L$ and by $L^*$ the dual of $L$.

\section{$K3$ surfaces associated with 8 points on the projective line}\label{k3}

\subsection{}\label{}
In this section, we shall construct a $K3$
surface associated to distinct 8 points from $P_1^8$.   In section \ref{discriminant locus}, 
we shall generalize this to the cases of any stable points and semi-stable points from $P_1^8$.
Let $\{ (\lambda_{i} : 1)\}$ be a set of distinct 8 points on the projective line.
Let $(x_{0}: x_{1}, y_{0}: y_{1})$ be the  bi-homogenious coordinates on ${\bf P}^{1} \times {\bf P}^{1}$.  
Consider a smooth divisor $C$ in ${\bf P}^{1} \times {\bf P}^{1}$ of bidegree $(4,2)$ given by
\begin{equation}\label{embed}
y_{0}^{2} \cdot \prod^{4}_{i=1} (x_{0} - 
\lambda_{i} x_{1}) + y_{1}^{2} \cdot 
\prod^{8}_{i=5} (x_{0} - \lambda_{i} x_{1}) = 0.
\end{equation}
Let $L_{0}$ (resp. $L_{1}$) be the divisor defined
by $y_{0} = 0$ (resp. $y_{1} = 0$).
Let $\iota$ be an involution of 
${\bf P}^{1} \times {\bf P}^{1}$ given by
\begin{equation}\label{invol}
(x_{0}: x_{1}, y_{0}: y_{1})  \longrightarrow
(x_{0}: x_{1}, y_{0}: -y_{1})
\end{equation}
which preserves $C$ and $L_{0}, L_{1}$.  
Note that the double cover of ${\bf P}^{1} \times {\bf P}^{1}$ 
branched along $C + L_{0} +L_{1}$ has 8 rational double points of type
$A_{1}$ and its minimal resolution $X$ is a $K3$ surface.   This $K3$ surface $X$ is obtained as follows:
First blow up the 8 points which are the intersection of $C$ and $L_0 + L_1$.  
Then $X$ is the double cover
branched along the proper transforms of $C$, $L_0$ and $L_1$.  
We remark that the isomorphism class of $X$ depends only on the 8 points in ${\bf P}^{1}$
 (i.e. independent on the order of 8 points) because elementary transformations change the order of 8 points.

The involution $\iota$ lifts to an automorphism $\sigma$ of order 4.
We can easily see that $\sigma^* \omega_X = \pm \sqrt{-1} \omega_X$ where $\omega_X$ is a
nowhere vanishing holomorphic 2-form on $X$.
We denote by $S_{0}, S_{1}$ the inverse image of
$L_{0}, L_{1}$ respectively.  
The projection
$$(x_{0}: x_{1}, y_{0}: y_{1})  \longrightarrow (x_{0}: x_{1})$$
from ${\bf P}^{1} \times {\bf P}^{1}$ to ${\bf P}^{1}$ induces an
elliptic fibration 
$$\pi : X \longrightarrow {\bf P}^{1}$$
which has 8 singular fibers of type $III$ in the sense of Kodaira  \cite{Kod} and
two sections $S_{0}, S_{1}$.
Let $E_{i} + F_{i}$ ($1 \leq i \leq 8$) be the
8 singular fibers of $\pi$.  Then we may assume that  
$$E_{i} \cdot S_{0} = F_{i} \cdot S_{1} = 1.$$
Put
\begin{equation}\label{}
H^{2}(X, {\bf Z})^+ = \{ x \in H^{2}(X, {\bf Z}) \mid (\sigma^{2})^{*}(x) = x \};
\end{equation}
\begin{equation}\label{}
H^{2}(X, {\bf Z})^- = \{ x \in H^{2}(X, {\bf Z}) \mid (\sigma^{2})^{*}(x) = -x \}.
\end{equation}
\noindent
We also denote by $S_X$, $T_X$ the {\it Picard lattice}, the {\it transcendental lattice} of $X$
respectively.  

\subsection{Lemma}\label{cohom} 
(i)   $H^{2}(X, {\bf Z})^+ \simeq U(2) \oplus D_{4} \oplus D_{4}$,  
$H^{2}(X, {\bf Z})^- \simeq U \oplus U(2) \oplus D_{4} \oplus D_{4}.$
\smallskip

(ii)  $\sigma^{*}$ {\it acts on $H^{2}(X, {\bf Z})^+$ trivially}. 
\smallskip

(iii)  {\it The following elements generate} $(H^{2}(X, {\bf Z})^+)^{*}/H^{2}(X, {\bf Z})^+ \simeq 
({\bf Z}/2{\bf Z})^{6}$:
$$(F_{1} + F_{2})/2, \
(F_{1} + F_{3})/2, \ (F_{1} + F_{4})/2, \ (F_{1} + F_{5})/2, \ 
(F_{1} + F_{6})/2, \ (F_{1} + F_{7})/2.$$

(iv)  {\it Let $U \oplus A_1^{\oplus 8}$ be the sublattice of $H^{2}(X, {\bf Z})^+$ generated by the
classes of a fiber, $S_0$ and $F_i$ $(1\leq i \leq 8)$.  Then $H^{2}(X, {\bf Z})^+$ is obtained from 
$U \oplus A_1^{\oplus 8}$ by adding the vector $(F_1 + \cdot \cdot \cdot + F_8)/2$.}

\begin{proof}
For the proof of the assertions (i)--(iii), see \cite{Kon1}, Lemma 5.2.   The sublattice 
$U \oplus A_1^{\oplus 8}$ has index 2 in $H^{2}(X, {\bf Z})^+$.  The later one is obtained from
the former by adding the class of $S_1$.  Hence the last assertion follows
from the fact that $S_1 = 2F + S_0 + (F_1 + \cdot \cdot \cdot + F_8)/2$ where $F$ is a fiber.
\end{proof}
\noindent
It follows that $H^{2}(X, {\bf Z})^+ \subset S_X$ and $T_X \subset H^{2}(X, {\bf Z})^-$.

\subsection{A quadratic form}\label{} 
First of all, we define: 
\begin{equation}
L=U^3\oplus E_8^3, \ M=U(2) \oplus D_{4} \oplus D_{4}, \  N=U \oplus U(2) \oplus D_{4}\oplus D_{4}.
\end{equation}

\noindent
Recall that $H^2(X, {\bf Z}) \cong L$.  We consider $M$ as a sublattice of $L$ and $N$ is the orthogonal
complement of $M$ in $L$.  It follows from Theorem 1.14.4 in Nikulin \cite{N3} that the embedding of $M$ into $L$ is unique.
Let $A_N = N^*/N$ which is isomorphic to a vector space ${\bf F}_{2}^6$ of dimension 6 over
${\bf F}_{2}$ (Lemma \ref{cohom}).
The {\it discriminant quadratic form }
\begin{equation}\label{quadratic}
q_{N} : A_{N}  \longrightarrow {\bf Q}/2{\bf Z}
\end{equation}
is defined by $q_N(x) = \langle x, x \rangle \ {\rm mod} \ 2{\bf Z}$.   In our situation, 
the image of $q_N$ is contained in ${\bf Z}/2{\bf Z}$ and hence $q_N$ is
a quadratic form  on $A_{N}$ defined over ${\bf F}_2$
whose associated bilinear form is given by
$$b_N(x,y) = 2 \langle x, y \rangle \ {\rm mod} \ 2{\bf Z}.$$
Let $u$ be the hyperbolic plane defined over ${\bf F}_2$, that is, 
the quadratic form of dimension 2 defined over ${\bf F}_2$ corresponding to the matrix
$
\begin{pmatrix}0&1
\\1&0
\end{pmatrix}
$.
The quadratic form $q_N$ is isomorphic to the direct sum of 3 copies of $u$: 
$q_N \cong u^{\oplus 3}$.  It is known that $(A_M=M^*/M, q_M) \cong (A_N, -q_N)$
(\cite{N3}, Corollary 1.6.2).

\subsection{Lemma}\label{} 
(i)   $O(q_{M}) \simeq O(q_{N}) \simeq O(u_{2}^{\oplus 3})
\simeq S_{8}$ {\it where $S_8$ is the symmetric group of degree $8$.}

(ii)  {\it The group $O(q_M)$ naturally isomorphic to the subgroup of $O(M)$ generated by
the permutations of the $8$ components $A_1$ in $U \oplus A_1^{\oplus 8}$.}

\begin{proof} 
The assertion (i) is well known (e.g. \cite{Atlas}, page 22).
For (ii), note that the permutations of 8 components $A_1$ can be extended to isometries of $M$
because they preserve  $(F_1 + \cdot \cdot \cdot + F_8)/2$ (see Lemma \ref{cohom}, (iv)).
Now the assertion is obvious.
\end{proof}

\subsection{A fundamental domain}\label{kaehler} 

Let $M = U(2) \oplus D_4 \oplus D_4$. 
Let $P(M)$ be a connected component of $\{ x \in M\otimes {\bf R} :  \ \langle x, x \rangle > 0 \}$.
Let $W(M)$ be the reflection group generated by $(-2)$-reflections
$$s_r : x \to x + \langle x, r \rangle r $$
for any $r \in M$ with $r^2 = -2$.  The group $W(M)$ acts on $P(M)$ discretely.  Let $C(M)$ be 
the finite polyhedral cone defined by the 18 $(-2)$-vectors which are corresponding to
the 18 smooth rational curves $S_0, S_1, E_i, F_i$ $(1 \leq i \leq 8)$ on $X$ under an isomorphism
$M \cong H^2(X, {\bf Z})^+$.

\subsection{Proposition}\label{fundamental domain} 
(i)  {\it The group $W(M)$ is of finite index in the orthogonal group $O(M)$.  Moreover the closure $\bar{C}(M)$ of $C(M)$ 
is a fundamental domain of $W(M)$.   The symmetry group of $C(M)$ is
isomorphic to $S_8\times {\bf Z}/2{\bf Z}$.} 

(ii) {\it If $S_X \cong M$, then $X$ contains exactly $18$ smooth rational
curves $S_0, S_1, E_i, F_i$ $(1 \leq i \leq 8)$.}

\begin{proof}
(i)  The first assertion follows from Nikulin's classification of such hyperbolic 2-elementary lattices (\cite{N1}, Theorem 4.4.1).   Moreover $C(M)$ satisfies
the condition in Vinberg's theorem \cite{V}, Theorem 2.6, that is, any maximal extended Dynkin diagram
in these 18 $(-2)$-vectors is either ${\tilde A}_1^{\oplus 8}$ or ${\tilde D}_4^{\oplus 2}$ both of which have the maximal rank 8.  
Hence $C(M)$  is of finite volume.  Let $C(M)'$ be a fundamental domain
of $W(M)$ with $C(M)'  \subset C(M)$.  Then \cite{V}, Lemma 2.4 implies that $C(M) = C(M)'$.
The last assertion is obvious.

(ii)  It follows from a remark in Vinberg \cite{V}, p. 335 that $C(M)$ has finite volume iff the polyherdal cone $\bar{C}(M)$ is
contained in the closure $\bar{P}(M)$ of $P(M)$.  If there exists a smooth rational curve $E$ different from 
the above 18 curves.
Then the intersection number of $E$ and any one of these 18 curves is non negative, that is, the class of $E$ is contained
in $\bar{C}(M) \subset \bar{P}(M)$.  This implies that $E^2 \geq 0$, which is a contradiction.
\end{proof}

\subsection{}\label{}
Let $X$ be a $K3$ surface as above.
Let $P(X)$ be the component of $\{ x \in S_X \otimes {\bf R} :  \ \langle x, x \rangle > 0 \}$
which contains an ample class.
Let $ \Delta (X)$ be the set of all effective divisors $r$ with $r^2 = -2$.
Let $C(X)$ be the polyhedral cone defined by:  
$$C(X) = \{ x \in S_X \otimes {\bf R} :  \ \langle x, r \rangle > 0 , r \in \Delta (X) \}.$$
Note that the integral points in $C(X)$ are nothing but the ample classes.  If $S_X \cong M$, then
$C(X) = C(M)$.

\subsection{Lemma}\label{absence} 
{\it The orthogonal complement of $H^2(X, {\bf Z})^+$ in $S_X$ contains no
$(-2)$-vectors.}

\begin{proof}
If $r\in (H^2(X, {\bf Z})^+)^{\perp} \cap S_X$ with $r^2 = -2$, then 
$(\sigma^*)^2(r) = -r$.  On the other hand, Riemann-Roch theorem implies that $r$ is effective.
This is a contradiction.
\end{proof}

\subsection{Proposition}\label{nodal} 
{\it Assume that $S_X = H^{2}(X, {\bf Z})^+$.  Then the automorphism group of $X$ is finite.  
Moreover $X$ has exactly $18$ smooth rational curves which are components of singular fibers of $\pi$ 
and two sections.}

\begin{proof}
Recall that $Aut(X)$ is isomorphic to $O(S_X)/W(S_X)$ up to finite groups (\cite{PS}).
Here $W(S_X)$ is the subgroup generated by $(-2)$-reflections.  Hence the assertion follows from
Proposition \ref{fundamental domain}.
\end{proof}

\subsection{The automorphism of order 4}\label{sigma} 
We shall study the action of $\sigma$ on $H^2(X, {\bf Z})^-$.  Recall that
$$D_4 \cong \{ (x_1,x_2,x_3,x_4) \in {\bf Z}^4 \ \mid \ x_1+x_2+x_3+x_4 \equiv 0 \ (\rm{mod}\ 2) \}.$$
Here we consider the standard inner product on ${\bf Z}^4$ with the negative sign.
Let $\rho_0$ be the isometry of $D_4$  given by
$$\rho_0(x_1,x_2,x_3,x_4) = (x_2,-x_1,x_4,-x_3).$$
Obviously $\rho_0$ is of order 4 and fixes no non-zero vectors in $D_4$.  
Also an easy calculation shows
that $\rho_0$ acts trivially on $D_4^*/D_4$.  Next let $e,f$ (resp. $e', f'$) be a basis of
$U$ (resp. $U(2)$).   Define the isometry $\rho_1$ of $U \oplus U(2)$ by
$$\rho_1(e) = -e-e', \ \rho_1(f) = f-f',  \ \rho_1(e') = e'+2e, \ \rho_1(f') = 2f-f'.$$
Obviously $\rho_1$ is of order 4, fixes no non-zero vectors in $U \oplus U(2)$ and acts trivially on
the discriminant group of $U \oplus U(2)$.  Thus we have an isometry
$\rho = \rho_1 \oplus \rho_0 \oplus \rho_0$ of $N = U\oplus U(2) \oplus D_4\oplus D_4$
which fixes no non-zero vectors in $N$ and acts trivially on $N^*/N$.  Then $\rho$ can be
extended to an isometry of the $K3$ lattice $L$ acting trivially on $M$ (\cite{N3}, Proposition 1.6.1).

\subsection{Lemma}\label{rho}  
{\it The isometry $\rho$ is conjugate to $\sigma^*$ under an isomorphism $H^2(X, {\bf Z}) \cong L$.}

\begin{proof}
Let $\omega$ be an eigenvector of $\rho$ which is sufficiently general, that is, satisfying the condition $\omega^{\perp} \cap L = M$.
By the surjectivity of the period map for $K3$ surfaces, 
there exists a $K3$ surface $Y$ and an isometry
$$\alpha_Y : H^2(Y, {\bf Z}) \to L$$
with $\alpha_Y (\omega_Y) = \omega$ where $\omega_Y$ is a nowhere vanishing holomorphic 
2-form on $Y$.   By the condition $\omega^{\perp} \cap L = M$, we have $S_Y \cong M$.
Consider the isometry $\phi =\alpha_Y^{-1} \circ \rho \circ \alpha_Y$ of $H^2(Y, {\bf Z})$.
Since $\rho$ acts trivially on $M$, $\phi$ preserves ample classes.  Then it follows from the Torelli
theorem \cite{PS} that $\phi$ is induced from an automorphism $g$ of $Y$ of order 4.
On the other hand, Proposition \ref{nodal} implies that 
$Y$ contains exactly 18 smooth rational curves
whose dual graph is the same as that of the smooth rational curves on $X$.  In particular,
$Y$ has an elliptic fibration with two sections and 8 singular fibers each of which is
type $III$ or $I_2$.  Since $g$ acts trivially on the Picard lattice $M$, $g$ preserves the
elliptic fibration and the class of each component of singular fibers of type $III$ or $I_2$.
Since the elliptic fibration has 8 singular fibers, $g$ acts trivially
on the base of the elliptic fibration, and hence induces an automorphism
of each fiber.  Hence all singular fibers are of type $III$.  
It follows from Nikulin \cite{N1}, Theorem 4.2.2 that the set of fixed points of the involution $g^2$ 
is the disjoint union of two smooth rational curves $R_0, R_1$ and a
smooth curve $C$ of genus 3.  Since $g$ acts trivially on the base, $R_0$, $R_1$ 
are sections of the
elliptic fibration.  We can easily see that $C$ passes through singular points of singular fibers of type $III$.
Thus we have the same configuration of
smooth rational curves on $Y$ as that of $X$.   By taking the quotient of $Y$ by $g^2$, we can see
that $Y$ is a deformation of $X$.  Hence we have the assertion.  
\end{proof}



\subsection{Markings}\label{marking} 

Recall that $H^2(X, {\bf Z})^+ \cong M = U(2) \oplus D_4 \oplus D_4$ (Lemma \ref{cohom}).  
We fix a fundamental domain $C(M)$  (Proposition \ref{fundamental domain}).  
It follows from Lemma \ref{rho} that there exists an isometry
 $$\alpha_X : H^2(X, {\bf Z}) \to L$$ 
satisfying $\rho = \alpha_X \circ \sigma^* \circ \alpha_X^{-1}$.
We call $\alpha_X$ a {\it marking} and the pair $(X, \alpha_X)$ a {\it marked} $K3$ surface.  
Then 

\subsection{Proposition}\label{ample cone} 
{\it There exists a marking $\alpha_X$ such that $\alpha_X (C(X)) \cap M\otimes {\bf R} 
\subset C(M)$.}

\begin{proof}
It follows from Lemma \ref{absence} that $\alpha_X (C(X)) \cap M\otimes {\bf R}$ is an open
polyhedral cone in $M\otimes {\bf R}$.
Hence Proposition \ref{fundamental domain} implies the assertion.
\end{proof}

\section{A complex ball uniformization}\label{unif}
In this section we construct an $S_8$-equivariant isomorphism between 
the moduli space of the projective equivalence classes of the set of distinct 
8 ordered  points in ${\bf P}^1$ 
and an open set of the arithmetic quotient of 5-dimensional complex ball.

\subsection{The period domain}\label{} 

Let $(X, \alpha_X)$ be a marked $K3$ surface and let
$\omega_X$ be a nowhere vanishing holomorphic 2-form on $X$.  Then
$\alpha_X(\omega_X) $ is contained in the following domain:

\begin{equation}\label{domain}
\calD = \{ \omega \in {\bf P}(N\otimes {\bf C}) : \langle \omega, \omega  \rangle =
0, \langle \omega, {\bar \omega} \rangle > 0 \}.
\end{equation}

\noindent
Note that $\calD$ is a disjoint union of two copies of a bounded symmetric domain of type $IV$ and of
dimension 10.  To get the period domain, we first define:

\begin{equation}\label{eigenspace}
V_{\pm} = \{ z \in N \otimes {\bf C} \mid \rho(z) = \pm \sqrt{-1} z \}.
\end{equation}

\noindent
It follows from Nikulin \cite{N2}, Theorem 3.1 that
$N \otimes {\bf C} = V_{+} \oplus V_{-}$.
Now we may assume $\sigma^*(\omega_X) = \sqrt{-1} \cdot \omega_X$.  Then
$\alpha_X(\omega_X) $ is, in fact, contained in $\calB$ defined by

\begin{equation}\label{domain2}
\calB = \{ z \in {\bf P}(V_{+}) \mid \langle z, \bar{z} \rangle > 0 \}.
\end{equation}

\noindent
If $z \in \calB$, then 
$$\langle z, z \rangle = \langle \rho(z),
\rho(z) \rangle = \langle \sqrt{-1}z, \sqrt{-1} z \rangle = -\langle
z, z \rangle,$$
and hence $\langle z, z \rangle = 0.$  Thus we have
$$\calD \cap {\bf P}(V_+) = \calB.$$
We remark that $\calB$ is a 5-dimensional complex ball.
We call $\alpha_X(\omega_X) $ the {\it period} of $(X, \alpha_X)$.
We also define two arithmetic subgroups:

\begin{equation}\label{group}
 \Gamma = \{ \gamma \in O(N) \mid \gamma \circ \rho = \rho \circ \gamma \};
\end{equation}

\begin{equation}\label{group}
 \Gamma(1-i) = {\rm Ker} ( \Gamma \to O(q_N) ).
\end{equation}

\smallskip
\noindent
We shall see that the quotient $\calB/\Gamma$ (resp. $\calB/\Gamma(1-i)$)
is the coarse moduli space of distinct 8 unordered points on ${\bf P}^1$
(resp. distinct 8 ordered  points on ${\bf P}^1$) (see Theorem \ref{main}).

\subsection{Hermitian form}\label{hermite}
We consider $N$ as a free ${\bf Z}[\sqrt{-1}]$-module $\Lambda$ by
$$(a+b\sqrt{-1})x = ax + b\rho(x).$$
Let 
$$h(x,y) = \sqrt{-1}\langle x, \rho(y) \rangle + \langle x, y \rangle.$$
Then $h(x,y)$ is a hermitian form on ${\bf Z}[\sqrt{-1}]$-module $\Lambda$.
With respect to a ${\bf Z}[\sqrt{-1}]$-basis $(1,-1,0,0)$, $(0,1,-1,0)$ of $D_4$,
the hermitian matrix of $h\mid D_4$ is given by

\begin{equation}\label{}
\begin{pmatrix}-2&1-\sqrt{-1}
\\1+\sqrt{-1}&-2
\end{pmatrix}.
\end{equation}

\noindent
And with respect to a ${\bf Z}[\sqrt{-1}]$-basis $e,e'$ of $U\oplus U(2)$,
the hermitian matrix of $h\mid U\oplus U(2)$ is given by

\begin{equation}\label{}
\begin{pmatrix}0&1+\sqrt{-1}
\\1-\sqrt{-1}&0
\end{pmatrix}.
\end{equation}

\noindent
Let 
$$\varphi : \Lambda \to N^*$$
be a linear map defined by $\varphi(x) = (x + \rho(x))/2.$
Note that $\varphi((1-\sqrt{-1})x) = \varphi(x - \rho(x)) = x \in N$. Hence
$\varphi$ induces an isomorphism
\begin{equation}\label{}
\Lambda/(1-\sqrt{-1})\Lambda \simeq N^*/N.
\end{equation}

\subsection{Remark}\label{M-Y} 
The hermitian form $h$ coincides with the one of Matsumoto and Yoshida in 
\cite{MY}, \S 6.  This implies that our groups $ \Gamma$, $ \Gamma(1-i)$ coincide
with the ones of Matsumoto and Yoshida in \cite{MY}.

\subsection{Reflections}\label{}
For $r \in N$ with $\langle r, r \rangle = -2$, we define a {\it reflection} 
$$s_r(x) = x + \langle r, x \rangle r$$
which is contained in $\tilde{O}(N) = {\rm Ker}(N \to O(q_N))$, but not in $\Gamma$.
On the other hand, by  considering $r$ as in $\Lambda$, we define a {\it reflection}

\begin{equation}\label{}
R_{r,\epsilon}(x) = x -(1-\epsilon){h(r,x)\over h(r,r)}r
\end{equation}
where $\epsilon \not= 1$ is a 4-th root of unity.  
We can easily see that $R_{r,-1}$ corresponds to the isometry in $\Gamma$
$$x \to x + \langle r, x \rangle r + \langle \rho(r), x \rangle \rho(r)$$
which coincides with $s_r \circ s_{\rho(r)}$.  Also
$R_{r,\sqrt{-1}}$ corresponds to the isometry in $\Gamma$
$$x \to x + \langle r, x \rangle (r-\rho(r))/2 + \langle \rho(r), x \rangle (r + \rho(r))/2$$
which induces a transvection of $A_N$ defined by
$$t_{\alpha}(x) = x + b_N(x,\alpha ) \alpha$$
where $\alpha \in A_N$ is a non-isotropic vector $ (r + \rho(r))/2 \ {\rm mod} \ N$.

%

\subsection{Discriminant}\label{discri} 

Let $r \in N$ with $r^{2} = -2$.
We denote by
$H_{r}$ the hyperplane of $\calB$ defined by
$$H_{r} = \{ z \in {\calB} : \langle z, r \rangle = 0 \}.$$
Let  ${\calH}$  be the union of all hyperplanes $H_{r}$
where $r$ moves on the set of all $(-2)$-vectors in $N$.  
We call $\calH$ the {\it discriminant locus}.
By Lemma \ref{absence}, the periods of marked $K3$ surfaces 
as above are contained in $\calB \setminus \calH$.  

Conversely let $\omega \in \calB \setminus \calH$.  Then by the surjectivity of the
period map, there exists a marked $K3$ surface $(X, \alpha_X)$ with
$\alpha_X (\omega_X) =  \omega$.   The condition $\omega  \notin \calH$ implies
that Proposition \ref{ample cone} holds for this $K3$ surface.  Hence, if necessary by replacing 
$\alpha_X$, we may assume that the isometry 
$\alpha_X^{-1} \circ \rho \circ \alpha_X$ preserves the ample cone of $X$.  It now follows
from the Torelli type theorem (\cite{PS})
that there exists an automorphism $\sigma$ of order 4
satisfying $\alpha_X^{-1} \circ \rho \circ \alpha_X = \sigma^*$.
Moreover  the marking defines an elliptic fibration 
$$\pi : X \to {\bf P}^1$$ with a section $s$.

\subsection{Lemma}\label{fibration}   
(i) {\it $\pi$ has $8$ singular fibers of type $III$ and two sections};

(ii) {\it The set of fixed points of $\sigma^2$ is the disjoint union of two sections and a
smooth curve of genus $3$ which passes through $8$ singular points of  $8$ singular fibers.}

\begin{proof}
It is known that the set of fixed points of the involution $\sigma^2$ is the disjoint union of
two smooth rational curves $R_0$, $R_1$ and a smooth curve $C$ of genus 3 
(Nikulin \cite{N1}, Theorem 4.2.2).  Obviously the set $X^{\sigma}$ of fixed points of  
$\sigma$ is contained in $R_0 + R_1 + C$.
Since $\sigma^* = \alpha_X^{-1} \circ \rho \circ \alpha_X$, $X^{\sigma}$ has the Euler number 12.  
Since $\sigma$ acts on $M$ trivially, it
preserves the section $s$ and the class of a fiber of $\pi$.  We show that $\sigma$ acts trivially
on the base of $\pi$.  Assume otherwise, then $X^{\sigma}$ is contained in two invariant
fibers, $F_1, F_2$.  Let $l$ be the number of irreducible one-dimensional components of
$X^{\sigma}$ and $k$ the number of isolated fixed points of $\sigma$.   Then
$2l + k = 12$.  
If we denote by $U$ the sublattice generated by the classes of a fiber and the section $s$, then
$U^{\perp} \cap  H^2(X, {\bf Z})^{\sigma^*} = A_1^8$.  Hence
the divisor $F_1 + F_2$ contains at least 10 components.  Assume $F_1$ contains at least
5 components.  Note that $\sigma$ preserves the component of $F_1$ which meets with $s$.
Obviously there are no singular fibers with non-trivial symmetry of order 4.  Hence the 
involution $\sigma^2$ preserves each component of $F_1$.  Then the sublattice
generated by components of $F_1$ not meeting $s$ has at least rank 4 and is isomorphic
to an indecomposable root lattice $R$.  Since $\sigma^2$ acts trivially on $R$, $R$ is contained
in $A_1^8$ which is impossible.

Thus $\sigma$ acts trivially on the base.  This implies that each fiber has an automorphism of
order 4.  In particular, singular fibers of $\pi$ are either of type $III, III^*$ or $I_0^*$.
Since
$U^{\perp} \cap  H^2(X, {\bf Z})^{\sigma^*} = A_1^8$ and
the singular fibers of type $III^*$ and $I_0^*$ have no non-trivial symmetry of order 4,
every singular fiber is of type $III$.  Since $\sigma$ fixes two points on each component of 
a singular fiber one of which is the singular point,  $R_0$, $R_1$ or $C$ passes through
these points.  Now we can easily see the assertion (ii).
\end{proof}

\subsection{Theorem}\label{main} 
 {\it The period map induces an $S_8$-equivariant  isomorphism $\phi$
between the moduli space $(P_1^8)^0$ of 
distinct ordered $8$ points on the projective line and the quotient space 
$(\calB \setminus \calH)/\Gamma(1-i)$.}

\begin{proof}  
As in \ref{discri}, for each $\omega \in  \calB \setminus \calH$, we have 
a marked $K3$ surface $(X, \alpha_X)$ with $\alpha_X(\omega_X) = \omega$.
Moreover $X$ has an automorphism $\sigma$ of order 4.
By Lemma \ref{fibration}, $X$ has an elliptic fibration with two section and 8 singular fibers of
type $III$.   By taking the quotient of $X$ by $\sigma^2$ and 
contracting $(-1)$-curves, we have the embedding of $C$ as in \ref{embed}.
This correspondence is the inverse of the period map.
\end{proof}

\section{Discriminant locus}\label{discriminant locus}

In this section we shall determine the discriminant locus $\calH$ of $\calB$. 

\subsection{}\label{} 
Let $r$ be a $(-2)$-vector in $N$.  Let $\omega \in \calB$.  Then
$\langle r, \omega \rangle = \sqrt{-1} \langle \rho(r), \omega \rangle$.  
Hence $r$ and $\rho(r)$ define
the same hyperplane $H_{r}$ in $\calB$.  Hence
$H_{r}$ corresponds to an embedding of the lattice
$R_r \cong A_{1} \oplus A_{1}$ generated by $r$ and $\rho(r)$ into $N$.
Obviously $R_{r} \simeq A_{1} \oplus A_{1}$.  Also every embedding of $R_{r}$ into $N$
is primitive, that is, $N/R_r$ is torsion free.

\subsection{Lemma}\label{A1}
 {\it Let $R_r^{\perp}$ be the orthogonal complement of $R_r$ in $N$.  Then 
$R_r^{\perp} \simeq U \oplus U(2) \oplus D_{4} \oplus A_{1}^{2}$.  In particular,
$(r + \rho(r))/2 \in N^*$.}

\begin{proof}
The proof for the first assertion is similar to those of  \cite{Kon1}, Lemmas 3.2, 3.3.
Then $R_r^{\perp} \oplus R_r$  is a sublattice of $N$ of index 2 and $N$ is obtained from
$R_r^{\perp} \oplus R_r$ by adding $(r + \rho(r))/2 + \theta$ where $\theta \in (R_r^{\perp})^*$.
We can see that $\langle (r+\rho(r))/2, x\rangle \in {\bf Z}$ for $x \in R_r \oplus R_r^{\perp}$ and
$x =  (r + \rho(r))/2 + \theta$.  Hence the second assertion holds.
\end{proof}

\subsection{}\label{64vectors}
Note that
$A_{N}$ consists of the following 64 vectors:
\smallskip

Type $(00): \alpha = 0, \ \# \alpha = 1$  ({\it zero});

Type $(0): \alpha \not=0, \ q(\alpha ) = 0, \ \# \alpha = 35$ ({\it non-zero isotropic vector});

Type $(1): q(\alpha ) = 1, \ \# \alpha = 28$ ({\it non-isotropic vector}).

\subsection{} \label{-4}
By Lemma \ref{A1}, the vector $(r + \rho(r))/2$ is contained in $N^*$.  In particular it defines a non-isotropic vector $(r + \rho(r))/2 \  {\rm mod} \ N$ in $A_N$.
Conversely let $\delta$ be a $(-4)$-vector in $N$ with $\delta/2 \in N^*$.
Since $\rho$ acts trivially on $A_N = N^*/N$, $\delta - \rho(\delta) \in 2N$.
Put $r = (\delta - \rho(\delta))/2 \in N$.  Since $\langle \delta, \rho(\delta) \rangle = 0$, 
$r^2 = -2$.  Obviously $\delta = r + \rho(r)$.  Thus we have

\subsection{Lemma}\label{-4vectors} 
$\{ \delta \in N \ : \ \delta^2 =-4, \ \delta/2 \in N^*\} = \{ r + \rho(r) \ :  \ r \in N, \ r^2 = -2 \}.$

\subsection{Lemma}\label{non-isot'} 
{\it Any non-isotropic vector in $A_N$ is represented by
$(r + \rho(r))/2$ for a suitable $(-2)$-vector $r$ in $N$.}

\begin{proof}
It follows from \cite{N3}, Theorem 1.14.2 that the natural map from $O(N)$ to $O(A_N)$ is
surjective.  The group $O(A_N) \cong S_8$
acts transitively on the set of non-isotropic vectors in $A_N$.
Combining these with Lemma \ref{-4vectors}, we have the assertion.
\end{proof}

\subsection{Proposition}\label{orbits}
(i)  $\Gamma/\Gamma(1-i)  \simeq S_8.$
\smallskip

(ii)  $\Gamma$ {\it acts transitively on the set of cusps of $\calB$ and on the set of
$\rho$-invariant  $R = A_1\oplus A_1$.}
\smallskip

(iii)   {\it $\Gamma(1-i)$-orbits of  $\rho$-invariant  $R = A_1\oplus A_1$ bijectively
correspond to non-isotropic vectors in $A_N$.  Also $\Gamma(1-i)$-orbits of cusps of
$\calB$ bijectively correspond to non-zero isotropic vectors in $A_N$.}

\begin{proof}
Recall that the pair $(N, \rho)$ naturally corresponds to the hermitian form $h$
(see Remark \ref{M-Y}).  Hence
the assertions follow from  \cite{MY}.  
\end{proof}

\subsection{Stable points}\label{}

Next we shall construct a $K3$ surface associated to each stable point from $P_1^8$.
Recall that 8 points is {\it stable} (resp. {\it semi-stable}) iff no four points (resp. five points) coincide
(e.g. \cite{DO}, Chap. I, \S 4, Example 2 (page 31)).
We denote by the symbol $(11111111)$ for distinct 8 points in ${\bf P}^1$.  If two points
(resp. three points) coincide, then we denote it by $(2111111)$ (resp. $(311111)$).

\subsection{Example: $(2111111)$}\label{}
We use the same notation as in  Section \ref{k3}.
We assume that $(\lambda_{1} : 1)$ has multiplicity 2.  Then the curve $C$ in \eqref{embed}
degenerates to one of the following two types:
$$C_1:  \ (x_{0} - \lambda_{1} x_{1})(y_{0}^{2} \prod^{4}_{i=2} (x_{0} - \lambda_{i} x_{1}) + 
y_{1}^{2} \cdot \prod^{7}_{i=5} (x_{0} - \lambda_{i} x_{1}) )= 0;$$
$$C_2: \  y_{0}^{2} (x_{0} - \lambda_{1} x_{1})^{2} \prod^{3}_{i=2} (x_{0} - 
\lambda_{i} x_{1}) + y_{1}^{2} \cdot 
\prod^{7}_{i=4} (x_{0} - \lambda_{i} x_{1}) = 0.$$
The minimal resolution $Y_i$ of the double covering of
${\bf P}^{1} \times {\bf P}^{1}$ branched along $C_i + L_{0} + L_{1}$
is a $K3$ surface and has an elliptic fibration $\pi$ which
has 6 singular fibers of type $III$, one singular fiber of type
$I_{0}^{*}$ and two sections $S_{0}, S_{1}$ ($i=1,2$).  We remark that $Y_1$ and $Y_2$ are
isomorphic because $C_1$ and $C_2$
are mutually transformed under elementary transformations.   Thus we denote by $Y$ instead of
$Y_1, Y_2$.  
Denote by $2R_{0} + R_{1} + R_{2} + R_{3} + R_{4}$ the singular fiber of
type $I_{0}^{*}$.  Assume that $S_{0}$ meets $R_{1}$ and $S_1$ meets $R_{2}$.
Then the normalization $\tilde C$ of $C_1$ meets $R_{3}, R_{4}$.  
The involution $\iota$ given in \eqref{invol} induces an automorphism $\sigma'$
of $Y$ of order 4.  Note that the restriction $\sigma'$ on 
$\tilde C$ is the hyperelliptic involution of the smooth
curve of genus 2.  This implies that $\sigma'$ switches $R_{3}$ and
$R_{4}$.  
Let $U$ be the sublattice generated by the classes of a fiber and $S_0$.  Then 
6 components in the fibers of type $III$ not meeting to $S_0$ and $R_2$, $2R_0 + R_2 + R_3 +
R_4$ generate the sublattice isomorphic to $A_1^8$.  This gives an isometry from $M$ into
$S_X$ such that $M^{\perp} \cap S_X$ contains $(-2)$ vectors $R_3$, $R_4$.  Thus the period of
$Y$ is contained in ${\calH}$.

For other stable points, the process is similar.  We have several types of the branch curve $C$ depending on the order of 8 points, however, they are transformed each other under 
elementary transformations.  Hence the corresponding $K3$ surface is determined by
the isomorphism class of 8 points (independent of the order of 8 points).

If three points coincide, then the corresponding elliptic fibration has a singular fiber of type $III^*$.
All cases except $(2222)$, the elliptic fibration has two sections.  In case of $(2222)$,
it has four sections.

The next Table 1 lists the type of 8 stable points on the projective line, type of singular fibers
of the elliptic fibration, the Picard lattice and the transcendental lattice of a generic member.

\begin{table}[h]
\[
\begin{array}{rllll}
{}& {\rm 8\ points}&{\rm Singular\ fibers}&{\rm Picard \ lattice}&{\rm Transcendental \ lattice} \\
1) & (1 1 1 1 1 1 1 1)&8 III&U(2)\oplus D_4 \oplus D_4& U\oplus U(2)  \oplus D_4\oplus D_4 \\
\noalign{\smallskip}
2) & (2 1 1 1 1 1 1) &I_0^*,\  6 III& U\oplus D_4\oplus D_4 \oplus A_1^{\oplus 2} &U\oplus U(2)  \oplus D_4\oplus A_1^{\oplus 2}  \\
    \noalign{\smallskip}
3) & (2 2 1 1 1 1)&2 I_0^*,\  4 III& U \oplus D_6 \oplus D_4 \oplus A_1^{\oplus 2} & U\oplus U(2)  \oplus  A_1^{\oplus 4}\\
    \noalign{\smallskip}
4) & (2 2 2 1 1) &3 I_0^*, 2\  III& U \oplus D_6 \oplus D_6 \oplus A_1^{\oplus 2} & A_1(-1)^{\oplus 2} \oplus  A_1^{\oplus 4} \\ 
\noalign{\smallskip}
5) & (2 2 2 2) &4 I_0^*& U \oplus D_8 \oplus D_8 & U(2) \oplus U(2)  \\ 
\noalign{\smallskip}
6) & (3 1 1 1 1 1) &III^*,\ 5 III& U \oplus D_8 \oplus D_4 & U\oplus U(2)  \oplus D_4\\
\noalign{\smallskip}
7) & (3 2 1 1 1) &III^*,\  I_0^*, \  III& U \oplus E_8 \oplus D_4 \oplus A_1^{\oplus 2} & 
U\oplus U(2)  \oplus A_1^{\oplus 2} \\ 
\noalign{\smallskip}
8) & (3 2 2 1) &III^*,\ 2 I_0^*, \ III& U \oplus E_8 \oplus D_6 \oplus A_1^{\oplus 2} &
A_1(-1)^{\oplus 2} \oplus  A_1^{\oplus 2}\\
\noalign{\smallskip}
9) & (3 3 1 1)&2 III^*, \ 2III& U \oplus E_8 \oplus  D_8 & U\oplus U(2)\\
    \noalign{\smallskip}
10) & (3 3 2) &2 III^*,\ I_0^*&  U \oplus E_8 \oplus D_{10} & A_1(-1)^{\oplus 2}\\
\end{array}
\]
\caption{}
\end{table}

\subsection{Strictly semi-stable points: $(44)$}\label{44}

In this case, we have the following 3 cases of curves in the quadric corresponding to
the strictly semi-stable points with unique minimal closed orbit:
$$C_3:  \ (x_{0} - \lambda_{1} x_{1})^2(x_{0} - \lambda_{2} x_{1})^2
(y_{0}^{2} + y_{1}^2)= 0;$$
$$C_4:  \ (x_{0} - \lambda_{1} x_{1})(x_{0} - \lambda_{2} x_{1})
(y_{0}^{2}(x_{0} - \lambda_{1} x_{1})^2 + 
y_{1}^2 (x_{0} - \lambda_{2} x_{1})^2)= 0;$$
$$C_5:  \ y_{0}^{2}(x_{0} - \lambda_{1} x_{1})^4
+ y_{1}^2 (x_{0} - \lambda_{2} x_{1})^4= 0.$$
These curves appear in the
list of Shah's classification of semistable $K3$ surfaces of degree 4.  See Shah \cite{S},
Theorem 4.8, B, Type II, (i)--(iii).

We denote by  $\bar{\calB}/\Gamma(1-i)$ the Satake-Baily-Borel compactification of 
${\calB}/\Gamma(1-i)$ whose boundary consists of 35 cusps.  Then we conclude:

\subsection{Theorem}\label{Main}
{\it The $S_8$-equivariant isomorphism $\phi$  in {\rm Theorem \ref{main}} can be extended to an 
$S_8$-equivariant
isomorphism $\tilde{\phi}$ between $P_1^8$ and $\bar{\calB}/\Gamma(1-i)$.  Moreover 
$\tilde{\phi}$ sends strictly semistable points to cusps and stable but not distinct $8$ points into
$\calH/ \Gamma(1-i)$}.

\begin{proof}
We can apply the argument of Horikawa's proof of the main theorem in \cite{H}.
Let $\calM$ be the space of all 8 semi-stable points on ${\bf P}^1$ and $\calM_0$ the space
of all distinct 8 points on ${\bf P}^1$.  We can easily see that $\calM \setminus \calM_0$ is locally
contained in a divisor with normal crossing.  By construction, $\phi$ is locally liftable to $\calB$.
It now follows from a theorem of Borel
\cite{Bo} that $\phi$ can be extended to a holomorphic map from $\calM$ to 
$\bar{\calB}/\Gamma(1-i)$ which induces a holomorphic map 
$$\tilde{\phi} : P_1^8 \to \bar{\calB}/\Gamma(1-i).$$
By using the same argument as in the proof of \cite{H}, Theorem 2.2, we can see that $\tilde{\phi}$ 
sends stable, but non-distinct 8 points to $\calH$.  More precisely we can choose a marking
for $K3$ surfaces corresponding to stable, but non-distinct 8 points, and define the period for
them.   For each stratification as in Table 1, we can prove the similar statement as in Lemma \ref{fibration}.
Then, 
as in the generic case (Theorem \ref{main}), by case by case argument according to strata,
we can see that the map $\tilde{\phi}$ is injective over $\calH$.
Moreover Shah's classification \cite{S} implies
the image of strictly semi-stable points go to the boundaries.
Hence the Zariski Main theorem implies that $\tilde{\phi}$ is an isomorphism.  
The $S_8$-equivariantness is obvious.
\end{proof}

\section{The Weil representation}\label{}

In this section we shall study the quadratic form $(A_{N}, q_{N})$
over ${\bf F}_{2}$ given in \eqref{quadratic} and the Weil representation of $SL(2, {\bf Z})$ on
the group ring ${\bf C}[A_{N}]$.

\subsection{}\label{}
Let

\begin{equation}\label{}
T =
\begin{pmatrix}1&1
\\0&1
\end{pmatrix},\quad
S =
\begin{pmatrix}0&-1
\\1&0
\end{pmatrix},
\end{equation}

\noindent
We denote by $\{ e_{\alpha} \}_{\alpha \in A_N}$ the satndard basis of ${\bf C}[A_N]$.
Let $\rho$ be the Weil representation of $SL(2,{\bf Z})$ on 
${\bf C}[A_{N}]$ which factors through $SL(2, {\bf Z}/2{\bf Z})$ (\cite{Bor}):

\begin{equation}\label{Weil}
\rho(T)(e_{\alpha}) = (-1)^{q_N(\alpha)} e_{\alpha};\quad
\rho(S)(e_{\alpha}) = {1 \over 8} \sum_{\beta} 
(-1)^{b_N(\beta, \alpha)} e_{\beta}.
\end{equation}

\noindent
Representatives of the conjugacy classes of $SL(2, {\bf Z}/2{\bf Z}) \simeq S_{3}$ consist 
of $ E, T, ST$.  A direct calculation shows that the traces of the action of $E, T, ST$
on ${\bf C}[A_{N}]$ are
$$tr(E) = 2^{6}, \ tr(T) = 8, \ tr(ST) = 1.$$
Let $\chi_{i}$ $(1 \leq i \leq 3)$ be the characters of irreducible
representations of $SL(2, {\bf Z}/2{\bf Z})$:
$\chi_{1}, \chi_{2}$ or $\chi_{3}$ is the trivial, alternating character or
the character of 2-dimensional irreducible representation respectively.
Let $\chi$ be the character of the Weil representation of
$SL(2, {\bf Z}/2{\bf Z})$ on ${\bf C}[A_{N}]$ and 
let $\chi = \sum_{i} m_{i} \chi_{i}$ be
its decomposition into irreducible characters.  Then
an elementary calculation shows that

\begin{equation}\label{decomposition}
\chi = 15 \chi_{1} + 7 \chi_{2} + 21 \chi_{3}.
\end{equation}
\smallskip
\noindent
We call a subspace $I$ of $A_{N}$ {\it totally isotropic} if
$q_{N}$ vanishes on $I$, and $I$ {\it maximal} if it has dimension 3.

\subsection{Lemma}\label{totally}
{\it For each maximal totally isotropic subspace 
$I$ of } $A_{N}$, 
$$\sum_{\alpha \in I} e_{\alpha} \in 
{\bf C}[A_{N}]^{SL(2, {\bf Z})}.$$

\begin{proof}
The proof is the same as that of  \cite{Kon2}, Lemma 3.2.  
\end{proof}

\smallskip
\noindent
Let $\alpha \in A_{N}$ with $q_{N}(\alpha) = 1$.  Then
$$t_{\alpha} : x \longrightarrow x + b_{N}(x, \alpha) \alpha$$
is called a {\it transvection} and contained in
$O(q_{N})$.  Note that $t_{\alpha}$ is
induced from a reflection $s_{r}$ associated with 
a $(-4)$-vector $r$ in $N$ 
with $r/2$ {\rm mod} $N = \alpha$ and these 
$t_{\alpha}$ ($\alpha \in A_{N}$ with $q_{N}(\alpha) = 1$) generate 
$O(q_{N})$.  These 28 transvections in $O(A_N)$ correspond to the 28 transpositions of $S_8$.

\subsection{Definition}\label{}
Let $q_s$ be a quadratic form on ${\bf F}_{2}^{3}$ given by
$$q_s(x) = \sum_{i=1}^{3} x_{i}, \ x = (x_{1},x_{2}, x_{3}) \in {\bf F}_{2}^{3}.$$
Note that the associated bilinear form of $q_s$ is identically zero.
Let $V$ be a 3-dimensional subspace of $A_{N}$.  We call $V$
{\it maximal totally singular} if $(V, q_{N} \mid V)$ is isomorphic to
$({\bf F}_{2}^{3}, q_s)$.  Obviously $V$ has a basis consisting of
3 mutually orthogonal non-isotropic vectors $\{ \alpha_{1}, \alpha_{2}, \alpha_{3} \}$.
We remark that $V$ consists of 4 non-isotropic vectors $\alpha_{i} \ (1\leq i \leq 3),
\alpha_{1} + \alpha_{2} + \alpha_{3}$ and
4 isotropic vectors $0, \alpha_{i} + \alpha_{j}, \ (1\leq i < j \leq 3)$. 
For each maximal totally singular subspace $V$, we define a vector $f_V \in {\bf C}[A_N]^{SL(2,{\bf Z})}$
on which the transvection $t_{\alpha}$ $(\alpha \in V)$ acts as $-1$.
Let $I$ be the kernel of $q_N \mid V$.  Then $I$ is a totally isotropic subspace of dimension 2 in $A_N$
and there exist exactly two maximal totally isotropic subspaces $I^+, I^-$ in $A_N$ which contain $I$.
We define

\begin{equation}\label{}
f_V = \sum_{\alpha \in I^+} e_{\alpha} - \sum_{\alpha\in I^-} e_{\alpha} \ \in {\bf C}[A_N].
\end{equation}

\subsection{Theorem}\label{unique}
{\it Let $V$ be a maximal totally singular
subspace.  Then $f_{V}$ is contained in
${\bf C}[A_{N}]^{SL(2, {\bf Z})}$ satisfying the following
condition}: 
$f_{V}$ {\it is the unique vector} ({\it up to constant}) {\it in}
${\bf C}[A_{N}]$ {\it on which transvections} $t_{\alpha}  (\alpha \in V, 
q_{N}(\alpha) = 1)$ {\it act as} $-1$.

\begin{proof}
The proof is the same as that of  \cite{Kon2}, Theorem 3.4.  
\end{proof}

\subsection{Remark-Definition}\label{W}
The group $O(q_{N}) (\simeq S_{8})$ naturally acts
on ${\bf C}[A_{N}]^{SL(2, {\bf Z})}$ with character
$\chi_{1} + \chi_{14}$, where $\chi_{1}$ is the trivial character and
$\chi_{14}$ is the character of an irreducible representation of
$S_{8}$ of degree 14.   This follows from Lemma \ref{totally} and \cite{Atlas}, page 22.
Moreover it follows from Theorem \ref{unique} that 
the multiplicity of the irreducible representation of degree
14 on ${\bf C}[A_{N}]$ is one.   We denote by $W$ the subspace of dimention
14 in  ${\bf C}[A_{N}]^{SL(2, {\bf Z})}$ with character
$\chi_{14}$.

\subsection{Lemma}\label{105}
 {\it The number of maximal totally singular 
subspaces of} $A_{N}$ {\it is equal to} $105$.

\begin{proof}
We can easily count the number of mutually orthogonal three non-isotropic vectors
in $A_N$ which is $2^3\cdot 3^2\cdot 5\cdot 7$.  On the other hand, the automorphism
group of a maximal totally singular subspace has order $2^3\cdot 3$.  Hence the assertion
follows. 
\end{proof}

\smallskip
\noindent
In the Lemma \ref{correspondence}, we shall give a geometric interpretation of maximal totally singular
subspaces.

\subsection{Heegner divisors}\label{}
Let $\delta \in N$ be a $(-4)$-vector with $\delta/2 \in N^*$.  
Let $D_{\delta}$ be the hyperplane of $\calD$ defined by
$$D_{\delta} = \delta^{\perp} \cap \calD.$$  
It follows from Lemma \ref{-4vectors} that
$$H_r = D_{\delta} \cap \calB$$
where $r = (\delta - \rho (\delta ))/2$ is a $(-2)$-vector in $N$ and $H_r$ is as in \ref{discri}.
For $\alpha \in A_N$ with $q_N(\alpha) = 1$, we define {\it Heegner divisors}
$\calD_{\alpha}$ and $\calH_{\alpha}$ by
$$\calD_{\alpha} = \sum_{\delta} D_{\delta}, \quad \calH_{\alpha} = \sum_{\delta} H_{\delta}$$
where $\delta$ varies over the set of $(-4)$-vectors in $N$ with $\delta/2 \ {\rm mod} \ N = \alpha$.
Since $H_r = D_{\delta} \cap \calB =  D_{\rho(\delta)} \cap \calB$, we have

\begin{equation}\label{heegner}
2\calH_{\alpha} =  \calD_{\alpha} \mid \calB.
\end{equation}

\section{Automorphic forms}\label{automorphic form}

\subsection{}\label{}
Let $\{ f_{\alpha} \}_{\alpha \in A_N}$ be a 
{\it vector valued elliptic modular form of weight $-4$ and of type}
$\rho$, i.e., $f_{\alpha}$ is a holomorphic function on the upper half plane satisfying

\begin{equation}\label{product}
f_{\alpha}(\tau + 1) = e^{2\pi i q(\alpha)} f_{\alpha}(\tau),
\quad
f_{\alpha}(-1/\tau) = {\tau^{-4}\over 8} \sum_{\beta \in A_N} e^{-2\pi i
\langle \alpha, \beta \rangle } f_{\beta}.
\end{equation}

\noindent
Recall that there are three types of vectors in $A_{N}$ denoted by
type 00, 0 or 1 according to zero, non-zero isotropic or non-isotropic respectively  
(see \ref{64vectors}) .
For each $\alpha \in A_{N}$, we denote by $m_{0}$ or $m_{1}$  the number 
of vectors $\beta \in A_{N}$ with $b(\alpha,\beta) = 0$ or $1$ 
respectively which is as in the following Table 2:

\begin{table}[h]
\[
\begin{array}{rllllllllllll}
\alpha& 00&00&00&0&0&0 &1&1&1\\
\beta&00&0&1&00&0&1&00&0&1\\
m_0&1&35&28&1&19&12&1&15&16\\
m_1&0&0&0&0&16&16&0&20&12\\
\end{array}
\]
\caption{}
\end{table}

\noindent
We shall find a modular form $h$ such that the components $h_{\alpha}$
are given by functions $h_{00}, h_{0}, h_{1}$ depending only on the
type of $\alpha$.  Then it follows from \eqref{product} and 
Table 2 that $h = \{ h_{\alpha} \}$ satisfies:

$$h_{00}(\tau + 1) = h_{00}(\tau) ,  \quad h_{00}(-1/\tau) = 
{\tau^{-4} \over 8}(h_{00}(\tau) + 35 h_{0}(\tau) + 28 h_{1}(\tau)),$$

$$h_{0}(\tau + 1) = h_{0}(\tau) ,  \quad h_{0}(-1/\tau) = 
{\tau^{-4} \over 8}(h_{00}(\tau) + 3 h_{0}(\tau) - 4 h_{1}(\tau)),$$

$$h_{1}(\tau + 1) = - h_{1}(\tau) , \quad  h_{1}(-1/\tau) = 
{\tau^{-4} \over 8}(h_{00}(\tau) - 5 h_{0}(\tau) + 4 h_{1}(\tau)).$$

\subsection{Lemma}\label{mod}
{\it One solution of these equations is given as
follows}:

$$h_{00}(\tau) = 56 \eta(2\tau)^{8}/\eta(\tau)^{16} =
56 + 896 q + 8064 q^{2} + \cdot \cdot \cdot ,$$ 

$$h_{0}(\tau) = -8 \eta(2\tau)^{8}/\eta(\tau)^{16} = -8 - 128 q 
- 1152 q^{2} - \cdot \cdot \cdot ,$$

$$h_{1}(\tau) = 8 \eta(2 \tau)^{8}/\eta(\tau)^{16} + 
\eta(\tau/2)^{8}/\eta(\tau)^{16} = q^{-1/2} + 36 q^{1/2} + 402 q^{3/2}
+ \cdot \cdot \cdot$$

{\it where} $\eta(\tau)$ {\it is the Dedekind eta function and}
$q = e^{2\pi \sqrt{-1}\tau}$.
\smallskip

\begin{proof}
The proof is the same as that of \cite{Kon2}, Lemma 4.3.
\end{proof}

\subsection{}\label{}
By applying Borcherds \cite{Bor},  Theorem 13.3, for the vector valued modular form $h$ given in 
Lemma \ref{mod}, we have

\subsection{Theorem}\label{Prod}
{\it There exists an automorphic form of weight $28 (=56/2)$
on} $\calD$ {\it which vanishes exactly on
Heegner divisors corresponding} 28 {\it non-isotoropic vectors in}
$A_{N}$.

\subsection{}\label{additive}
On the other hand, by Borcherds \cite{Bor},  Theorem 14.3, we have an
$S_8$-equivariant map 
$$\varphi : W  \to \calA_4(\tilde{O}(N))$$ 
where $W$ is the 14-dimensional subspace of ${\bf C}[A_{N}]^{SL(2, {\bf Z})}$ 
given in Remark-Definition \ref{W}, $\calA_4(\tilde{O}(N))$ is
the space of automorphic forms on $\calD$ of weight 4 with respect to $\tilde{O}(N)$ and
$$\tilde{O}(N) = {\rm Ker}(O(N) \to O(q_N)).$$  It follows from \cite{N3}, Theorem 1.14.2 that
the map from $O(N)$ to $O(q_N)$ is surjective and hence $S_8$ $(\cong O(q_N) \cong O(N)/\tilde{O}(N))$
naturally acts on $\calA_4(\tilde{O}(N))$.  On the other hand, $S_8$ acts on $W$.  With respect to
these actions, $\varphi$ is $S_8$-equivariant.

\subsection{Lemma}\label{} {\it The map $\varphi$ is injective}.

\begin{proof}
The proof is the same as that of \cite{Kon2}, Lemma 4.1.
\end{proof}

\subsection{Theorem}\label{lift}
 {\it Let} $V$ {\it be a maximal totally singular
subspace of} $A_{N}$.  {\it Let} $F_{V}$ {\it be an automorphic form
of weight $4$ on $\calD$
associated with} $f_{V} \in {\bf C}[A_{N}]^{SL(2, {\bf Z})} : \ F_V = \varphi (f_V)$.
{\it Then}
$$(F_{V}) = \sum_{\alpha \in V, \ q(\alpha) = 1} {\calD}_{\alpha}$$
{\it where ${\calD}_{\alpha}$ is the Heegner divisor associated with} $\alpha$.

\begin{proof}
Let $\Phi$ be the product of all $F_{V}$ where $V$ varies over all maximal totally singular 
subspaces.  Then $\Phi$ is an automorphic form of weight $105 \times 4$ (Lemma \ref{105}).  
By Theorem \ref{unique} and the $S_8$-equivariantness of $\phi$, we have
$$\sum_{\alpha \in V, \ q(\alpha) = 1} {\calD}_{\alpha} \subset (F_V).$$
Hence $\Phi$ vanishes along
${\calD}_{\alpha}$ with vanishing order $\geq 4 \times 105 / 28 = 15$.
On the other hand, the $15$-th power of the automorphic form given in
Theorem \ref{Prod} has the same weight and vanishes on ${\calD}_{\alpha}$ with
multiplicity 15.  The assertion now follows from the Koecher principle.
\end{proof}

\section{Cross ratios}\label{}

\subsection{}\label{}

Let

\begin{equation}\label{tau}
\tau =
\begin{pmatrix} \tau_{11} & \tau_{12}
\\\tau_{21} & \tau_{22}
\\\tau_{31} & \tau_{32}
\\\tau_{41} & \tau_{42}
\end{pmatrix},\quad
\tau_{ij} \in \{ 1,2,..., 8 \}.
\end{equation}
\noindent
We call $\tau$ a {\it tableau}.  A tableau $\tau$ is called {\it standard} if 
$$\tau_{ij} < \tau_{ij+1}, \quad \tau_{ij} \leq \tau_{i+1 j}$$
for any $i, j$.  The number of tableaus is 105 and the number of standard tableaus is 14
(e.g. \cite{DO}, Chap. I).
For each $\tau$ we define
$$\mu_{\tau} = \prod_{1 \leq i \leq 4} {\rm det}(v^{\tau_{i1}} v^{\tau_{i2}})$$
where $v^{i} \in {\bf C}^{2}$ is a column vector.  If

\begin{equation}\label{}
v^{i}  =
\begin{pmatrix} 1
\\x^{i} 
\end{pmatrix},
\end{equation}
\noindent
then, 
$$\mu_{\tau} = \prod_{1 \leq i \leq 4} (x^{\tau_{i2}} - x^{\tau_{i1}}).$$
These $\mu_{\tau}$ define an $S_8$-equivariant map $\Theta$ from 
$P_1^8$ to ${\bf P}^{13}$ (e.g. \cite{DO}).
We identify  $P_1^8$ with 
$\bar\calB/\Gamma(1-i)$ under the isomorphism given in Theorem \ref{Main}.  In the following,
we shall discuss a relation between $\Theta$ and the map defined by 14-dimensional space 
$W$ of automorphic forms given in \S \ref{automorphic form}.

\subsection{}\label{}

We give a relation between the set of tableaus and the set of totally singular subspaces.
Let $K= U \oplus A_1^{\oplus 8}$ be the sublattice given in Lemma \ref{cohom}.
Then $A_K = K^*/K \cong ({\bf F}_2)^8$ is generated by $F_i/2$ $(1\leq i \leq 8)$ corresponding to
8 points on the projective line.
The discriminant quadratic form $q_K$ of $K$ is a map
$$q_K : ( {\bf F}_2)^8 \to {\bf Q}/2{\bf Z}$$
defined by $q_K(x) = \langle x, x \rangle \ \rm{mod} \ 2{\bf Z}$.
Let $\theta = (F_1 +  \cdot  \cdot \cdot + F_8)/2$ which is perpendicular to all vectors in $A_K$.
Then it is known (Nikulin \cite{N3}, Proposition 1.4.1)
that the discriminant quadratic form of $M$ is obtained by
$$q_M = q_K \mid \theta^{\perp} / \theta.$$
Finally $q_N$ is canonically isomorphic to $-q_M$ (\cite{N3}, Corollary 1.6.2).
Thus non-isotropic vectors with respect to
$q_N$ bijectively corresponds to the vectors $(F_i + F_j)/2$ $(i \not= j)$ in $A_K$.
Each column $(\tau_{i1}, \tau_{i2})$ of $\tau$ in \eqref{tau} 
defines a vector in $({\bf F}_2)^8$  whose nonzero entries are indexed by $\tau_{i1}, \tau_{i2}$.
Thus four columns of $\tau$ corresponds to mutually orthogonal 4 non-isotropic vectors
which generate a maximal totally singular subspace in $A_N$.
This implies the following:

\subsection{Lemma}\label{correspondence}
{\it The set of 105 tableaus $\tau$ bijectively corresponds to the set of maximal totally singular
subspaces of $A_N$.  Under this correspondence, the zero of $\mu_{\tau}$ coincides with the zero of $F_{V}$ where $V$ is a maximal totally singular subspace corresponding to $\tau$.}

\subsection{}\label{}

Consider the linear system of automorphic forms $F_V$ of dimension 14 defined by $W$
(see \ref{additive}).
Note that the divisor $(F_V \mid \calB)$ is given by
$$2 \sum_{\alpha \in V, \ q(\alpha) = 1} {\calH}_{\alpha}$$
(see \eqref{heegner}, Theorem \ref{lift}).
Since $\calB$ is simply connected, we can take
a square root of $F_V$.  Thus we have an automorphic form $G_V$ of weight 2 on $\calB$ with

\begin{equation}\label{restriction}
(G_V) = \sum_{\alpha \in V, \ q(\alpha) = 1} {\calH}_{\alpha}.
\end{equation}
\noindent
Then $\{ G_V \}_V$ defines a map $\Psi$ from $\bar{\calB}/\Gamma(1-i)$ to 
${\bf P}^{13}$.

\subsection{Theorem}\label{cross}
{\it The map $\Theta$ coincides with $\Psi$.}

\begin{proof}
We shall show that 
$\mu_{\tau_1}/\mu_{\tau_2}$ coincides with $G_{V_1}/G_{V_2}$
for suitable tableaux $\tau_1, \tau_2$ and the corresponding 
maximal totally singular subspaces $V_1, V_2$.
We consider the following tableaux:

\begin{equation}\label{}
\tau_1 =
\begin{pmatrix} 1 & 2
\\3 & 4
\\5 & 6
\\7 & 8
\end{pmatrix},
\tau_2 =
\begin{pmatrix} 1 & 2
\\3 & 4
\\5 & 7
\\6& 8
\end{pmatrix},
\tau_3 =
\begin{pmatrix} 1 & 2
\\3 & 4
\\5 & 8
\\6 & 7
\end{pmatrix}.
\end{equation}
\noindent
We take a decomposition of $A_N = u_1\oplus u_2 \oplus u_3$ into
three hyperbolic planes $u_1, u_2, u_3$ defined over ${\bf F}_2$.  Let $\{e_i, f_i\}$ be a basis of $u_i$ with
$\langle e_i, e_i \rangle = 0, \ \langle f_i, f_i \rangle = 0, \ \langle e_i, f_i \rangle = 1.$
Let $\alpha_i = e_i + f_i$ be the non-isotropic vector in $u_i$.  
We may assume that
$$V_1 = \langle \alpha_1, \alpha_2, \alpha_3\rangle, \ 
V_2 = \langle \alpha_1, \alpha_2, \alpha_1 + e_3\rangle, \ 
V_3 = \langle \alpha_1, \alpha_2, \alpha_1+f_3\rangle$$
correspond to $\tau_1, \tau_2,  \tau_3$ respectively.  
Then by \eqref{restriction}, we can see that 
$(G_{V_1}/G_{V_2}) = (\mu_{\tau_1}/\mu_{\tau_2})$ as divisors.
Note that the function $\mu_{\tau_1}/\mu_{\tau_2}$ takes the value 1 on the divisors defined by $x_5-x_8$
and $x_6-x_7$.
On the other hand, easy calculation shows that $f_{V_1} - f_{V_2} = f_{V_3}$ and hence
$F_{V_1} - F_{V_2}$ vanishes on the divisor of $F_{V_3}$.  This implies that $G_{V_1} / G_{V_2} = 1$ on the divisors $\calH_{a_1+f_3}$ and $\calH_{a_2+f_3}$.  
Hence $G_{V_1}/G_{V_2} = \mu_{\tau_1}/\mu_{\tau_2}$.
For any pair of $\tau_1', \tau_2'$, the ratio $\mu_{\tau_1'}/\mu_{\tau_2'}$  can be
written as the product of some $\mu_{\tau_1}/\mu_{\tau_2}$ as above type.
Hence the assertion follows.
\end{proof}

\subsection{Theorem}\label{morphism}
{\it 
 $\Psi$ is an embedding from $\bar{\calB}/\Gamma(1-i)$ into ${\bf P}^{13}$.  The image satisfies $2^2\cdot 3\cdot 5 \cdot 7$ quartic relations.}

\begin{proof}
It is known that  $\Theta$ is embedding (Koike \cite{Koi}).  
The proof of the second assertion is the same as that of \cite{Kon2}, Theorem 7.2,
that is, for each non-isotropic vector in $A_N$ we have 15 quartic relations.
Since the number of non-isotropic vectors is 28, the second assertion follows.
\end{proof}

\subsection{Remark}\label{}  In the paper \cite{Koi}, by using a computer, he showed that the image is the intersection of 14 quadrics.

\subsection{Remark}
In the paper \cite{MT}, Matsumoto and Terasoma constructed an $S_8$-equivariant map from
the moduli space of ordered 8 points on ${\bf P}^1$ to ${\bf P}^{104}$ by using the theta
constants related to the curve which is the 4-fold covering of ${\bf P}^1$ branched at
8 points.  They showed that this map coincides with the map defined by the above 
105 $\mu_{\tau}$.

\subsection{Remark}
Since 8 points on ${\bf P}^1$ naturally correspond to hyperelliptic curves of genus three, we can consider 
that our case is a degenerate one of 
smooth curves of genus three.  The moduli space of non-hyperelliptic curves of genus three can be also
described as an arithmetic quotient of a complex ball (\cite{Kon1}).   On the other hand, Coble constructed
a map from the moduli space of curves of genus 3 with level 2-structure to ${\bf P}^{14}$ 
by using G{\" o}pel functions (Coble \cite{C}, Dolgachev, Ortland \cite{DO}, Chap. IX).
It would be interesting to extend the result in this paper to the case of curves of genus three.



\smallskip


\begin{thebibliography}{AMR}


\bibitem[AF]{AF} D.\ Allcock, E.\ Freitag, {\it
Cubic surfaces and Borcherds products},
Comm. Math. Helv., {\bf 77} (2002), 270--296.


\bibitem[Atlas]{Atlas} J.\ H.\ Conway et al., {\it Atlas of finite groups},
Oxford  1985.


\bibitem[Bor]{Bor} R.\ Borcherds,  {\it Automorphic forms with singularities on
Grassmannians}, Invent. Math.  {\bf 132} (1998), 491--562.

\bibitem[Bo]{Bo} A.\ Borel,  {\it Some metric properties of arithmetic quotients of symmetric
spaces and an extension theorem}, J. Diff. Geometry  {\bf 6} (1972), 543--560.

\bibitem[C]{C} A.\ Coble,  {\it Algebraic geometry and theta functions}, Amer. Math. Soc. Coll. Publ.  {\bf 10} Providence, R.I., 1929 (3rd ed., 1969).

\bibitem[DM]{DM} P.\ Deligne, G.\ W.\ Mostow, {\it Monodromy of
hypergeometric functions and non-lattice integral monodromy}, Publ. Math.
IHES, {\bf 63} (1986), 5--89

\bibitem[DO]{DO} I.\ Dolgachev, D. Ortland, {\it Point sets in projective spaces and
theta functions}, Ast{\'e}risque {\bf 165}(1988).

\bibitem[H]{H} E.\ Horikawa,  {\it On the periods of Enriques surfaces.} II, Math. Ann.  {\bf 235} (1978), 217--246.

\bibitem[Kod]{Kod}  K.\ Kodaira, {\it On compact complex analytic surfaces} II, Ann. Math.,
{\bf 77}(1963),  563--626.  III, Ann. Math., {\bf 78}(1963),  1--40.

\bibitem[Koi]{Koi}  K.\ Koike, {\it The projective embedding of the configuration space 
$X(2,8)$}, preprint.

\bibitem[Kon1]{Kon1}  S.\ Kond${\rm \bar o}$,
{\it  A complex hyperbolic structure for 
the moduli space of curves of genus three},  J. reine angew. Math., 
{\bf 525}(2000),  219--232.

\bibitem[Kon2]{Kon2} S.\ Kond${\rm \bar o}$, {\it The moduli space of Enriques surfaces
and Borcherds products},  J. Algebraic Geometry, {\bf 11} (2002), 601--627.


\bibitem[MT]{MT} K.\ Matsumoto, T.\ Terasoma, {\it Theta constants associated to coverings of
${\bf P}^1$ branching at $8$ points},  Compositio Math., {\bf 140} (2004), 1277--1301.

\bibitem[MY]{MY} K.\ Matsumoto, M.\ Yoshida, {\it Configuration space of $8$ points 
on the projective line and a $5$-dimensional Picard modular group}, 
Compositio Math., {\bf 86} (1993), 265--280.


\bibitem[N1]{N1} V.\ V.\ Nikulin, {\it Factor groups of groups of
automorphisms of hyperbolic forms with
respect to subgroups generated by $2$-reflections}, J. Soviet Math., {\bf
22} (1983), 1401--1475.

\bibitem[N2]{N2} V.\ V.\ Nikulin, {\it Finite automorphism groups of K{\" a}hler $K3$ 
surfaces}, Trans. Moscow Math. Soc., {\bf 38} (1980), 71--135.

\bibitem[N3]{N3} V.\ V.\ Nikulin, {\it Integral symmetric bilinear forms
and its applications}, Math. USSR Izv., {\bf 14} (1980), 103--167.

\bibitem[PS]{PS} I.\ Piatetski-Shapiro, I.\ R.\ Shafarevich, {\it A
Torelli theorem  for algebraic surfaces of type $K3$}, Math. USSR Izv.,
{\bf 5} (1971), 547--587.

\bibitem[S]{S} J.\  Shah, {\it Degenerations of $K3$ surfaces of degree $4$}, 
Trans. A. M. S., {\bf 263} (1981), 271--308.

\bibitem[V]{V} E.\ B.\ Vinberg, {\it Some arithmetic discrete groups in Lobachevskii spaces}, 
in "Discrete subgroups of Lie groups and applications to moduli", Tata-Oxford (1975), 323--348.

\end{thebibliography}
\end{document}